\documentstyle{bcp98e}
\input{amssym.def}
\input{amssym}

\renewcommand{\thesection}{\arabic{section}}
\newtheorem{thm}{\arabic{section}.\arabic{abz}. Theorem \hspace{-3mm}}

\newtheorem{defi}{\arabic{section}.\arabic{abz}. Definition \hspace{-3mm}}

\newtheorem{prop}{\arabic{section}.\arabic{abz}. Proposition \hspace{-3mm}}

\newtheorem{cor}{\arabic{section}.\arabic{abz}. Corollary \hspace{-3mm}}

\newtheorem{exa}{\arabic{section}.\arabic{abz}. Example \hspace{-3mm}}

\newtheorem{rem}{\arabic{section}.\arabic{abz}. Remark \hspace{-3mm}}

\newcounter{abz}[section]
\newcounter{sabz}[abz]
\addtocounter{section}{-1}
\newcommand{\abz}{\refstepcounter{abz}}
\newcommand{\sabz}{\refstepcounter{sabz}}

\newcommand{\lbr}{\linebreak[0]}
\def\pr{{\Bbb P}}
\def\d{\partial}
\def\Ci{C^{\infty}}
\def\G{\Gamma}
\def\la{\lambda}
\def\ra{\mathop{\rm rank}\nolimits}
\def\ke{\mathop{\rm ker}\nolimits}

\def\Span{\mathop{\rm Span}\nolimits}
\def\Sing{\mathop{\rm Sing}\nolimits}
\def\i{{\rm i}}
\def\codim{\mathop{\rm codim}\nolimits}
\def\R{{\Bbb R}}
\def\C{{\Bbb C}}
\def\A{{\frak g}}
\def\D{{\frak g}^*}
\def\V{{\cal V}}
\def\W{{\cal W}}
\def\S{{\cal S}}

\begin{document}
\mathclass{58F07,53A60}

\abbrevauthors{A.Panasyuk}
\abbrevtitle{Veronese webs of higher codimension}

\title{Veronese webs for bihamiltonian structures \\
of higher corank}
\author{Andriy\ Panasyuk}

\address{Division of Mathematical Methods
in Physics, \\ University of Warsaw, \\ Ho\.{z}a St. 74, 00-682
Warsaw, Poland\\
E-mail: panas@fuw.edu.pl\\ and \\ Mathematical Institute of
the Polish Academy of Sciences\\ \'{S}niadeckich St. 8, 00-950
Warsaw, Poland}

\vfootnote{}{Partially supported by the Polish grant KBN 2 PO3A 135 16.}
\maketitlebcp

\begin{tabbing}
\hspace{7.3cm}\= {\em To the memory of Stanis\l aw Zakrzewski,}\\
\> {\em with  the respect and gratitude}
\end{tabbing}

\section{Introduction.}

A $\Ci$- manifold $M$ is endowed by a Poisson
pair if two linearly independent smooth bivectors
$c_1,c_2$ are defined on $M$
and  $c_\la
=\la_1c_1+\la_2c_2$ is a  Poisson bivector for any $\la
=(\la_1,\la_2)\in\R^2$. A bihamiltonian structure
$J=\{c_\la\}$  is the whole 2-dimensional family of bivectors.
The structure $J$ is degenerate if $\ra
c_\la<\dim\ M, \la\in\R^2$.

An intensive study of such objects was done by
I.M.Gelfand and I.S.Zakharevich (\cite{gz1}, \cite{gz2},
\cite{gz3}) in a particular case of bihamiltonian structures in
general position on an odd-dimensional $M$ (the corresponding
Poisson pairs are necessarily degenerate: $\ra\, c_\la=2n, \la\in
\R^2\setminus\{0\}$, if $\dim M=2n+1$).
In \cite{gz2} there was introduced a notion of a Veronese web,
i.e. a 1-parameter family
of 1-codimensional foliations such that the corresponding family
of annihilators is represented by the Veronese curve in the
cotangent space at each point. It turns out that Veronese webs
form a complete system of local invariants for bihamiltonian structures
of general position. More precisely, it
was shown in \cite{gz2} that any such structure
$J=\{c_\la\}$ in $\R^{2n+1}$ admits a
local reduction to a Veronese web ${\cal W}_J$ on a
$(n+1)$-dimensional manifold and that for any Veronese web
${\cal W}$ one can locally construct a bihamiltonian
structure $J({\cal W})$ of general position in $\R^{2n+1}$
with the reduction equal to ${\cal W}$. In the real
analytic case $J$ and $J({\cal W}_J)$ are isomorphic.

The aim of this paper is to introduce a wider class of
degenerate bihamiltonian structures that possess many
features of the general position case and to generalize the
notion of a Veronese web for this class. We call the
bihamiltonian structures from this class complete since
they are intemately connected with the completely
integrable systems (\cite{arnold}) on $M$.  In particular,
the Poisson pairs appearing in the well known method of
argument translation (see \cite{mf}, \cite{fom}, and
Example \ref{a20.50}, below) generate complete
bihamiltonian structures of higher ($>1$) corank.

The paper is organized as follows. In Section 1 we recall some
definitions and facts about bihamiltonian structures and
introduce the main definition of completeness. The last is based
on one result of A.Brailov (Theorem \ref{a20.30}). We show that
complete bihamiltonian structures generalize the case of general
position. Analyzing the corresponding Poisson pair
$(c_1(x),c_2(x))$ at a point $x\in M$ we deduce that it consists
of finite number of the so called Kronecker blocks (Corollary
\ref{a20.67}); the general position is characterized by the case
of the sole block. Section 2 is devoted to distinguishing the
invariants for the sum of $k$ Kronecker blocks. In the next
section we define local Veronese webs for complete bihamiltonian
structures under some additional assumption of simplicity. This last
means that: 1)the number of Kronecker blocks does not change from
point to point and the corresponding subspaces vary smoothly
"sweeping" a flag of $k$ subbundles in the tangent bundle; 2)there are no
blocks of equal dimension. The second condition allows to avoid some
technical complications but in principle may be skiped
(see Remark~\ref{rema}).  In general, the mentioned
distributions are nonintegrable (Examples \ref{40.40},
\ref{40.50}); consequently, the bihamiltonian structure
does not split to direct product of the bihamiltonian
structures of corank $1$, i.e.  of general position. We
conclude the paper calculating the Veronese web for the
method of argument translation (Section 5).  In the case of
normal noncompact real form of complex simple Lie algebra
this web is generically a product of flat Veronese webs of
codimension $1$.

Recent papers \cite{gz4}, \cite{gz5} are closely related to
the subject, in particular to generalized Veronese webs.
In \cite{gz5} the author introduces a more general notion
of a Kronecker web, which is essentially equivalent to the
notion of a Veronese web (see Definition~\ref{40.20}) in
case of simple bihamiltonian structures.  Our approach
emphasizes a bit more the role of Veronese curves in the
theory.

The following two questions arise from the context of this paper.

1. Does the Veronese web of a complete bihamiltonian
structure determine it up to an isomorphism?

2. What is a  relation between the Veronese webs
introduced here and $d$-webs of maximal rank  and codimension $2$
studied in paper \cite{cg1} of S.S.Chern and P.A.Griffiths? (The
notion of the rank of a $d$-web should not be confused with
that of a bihamiltonian structure; corank of bihamiltonian
structure is equal to  codimension of the web.)

Note that the $d$-webs of maximal rank and codimension $1$
considered in  paper \cite{cg2} of the same authors are
intemately connected with the Veronese webs of codimension $1$.

The author would like to thank Prof.
Ilya Zakharevich for useful remarks on this paper and for
indicating references \cite{gz4},~\cite{gz5}, which had an essential inluence
on its final version.

\section{Bihamiltonian structures and completeness.}
Let $M$ be a $\Ci$- manifold.
In the sequel, all considered Poisson bivectors will
have
maximal rank on an open
dense subset in $M$.
Given a Poisson bivector $c$, define $\ra c$ as
$\max_{x\in M}\ra c(x)$.

\abz\label{a10.180}
\begin{defi} \rm
Two linearly independent  Poisson bivectors $c_1,c_2$
on $M$ form a {\em Poisson pair} if $c_\la=\la_1c_1+\la_2c_2$
is  a  Poisson bivector for any $\la=(\la_1,\la_2)\in\R^2$.
\end{defi}

\abz\label{a10.190}
\begin{prop}
A pair of linearly independent Poisson bivectors $(c_1,c_2)$ is
Poisson if and only if $[c_1,c_2]=0$, where $[\cdot,\cdot]$ is
the Schouten bracket.
\end{prop}

\abz\label{a10.210}
\begin{defi} \rm
A {\em bihamiltonian structure }
on $M$ is defined as a two-dimensional linear subspace
$J=\{c_\la\}_{\la\in{\cal S}}$ of Poisson bivectors on $M$
parametrized by a two-dimensional vector space ${\cal S}$ over
$\R$ .  We say that
$J$ is {\em degenerate} if $\ra c_\la<\dim M$ for any $c_\la\in J$.
\end{defi}

It is clear that every Poisson pair generates a bihamiltonian
structure and the transition from the latter one  to a Poisson
pair corresponds to a choice of basis in ${\cal S}$. We shall write
$(J,c_1,c_2)$ for a bihamiltonian structure $J$ with a chosen
Poisson pair $(c_1,c_2)$ generating $J$.

\abz\label{a10.230}
\begin{defi} \rm
Let $J$ be a bihamiltonian structure. Introduce a {\em subfamily}
$J_0\subset J$ of Poisson bivectors of
maximal rank $R_0$ (the set $J\setminus J_0$ is at most a finite sum of
1-dimensional subspaces), and a {\em set} of functions ${\cal
F}_0=\Span_\R (\bigcup_{c\in J_0}Z_c(M))$, where $Z_c(M)$
stands for the space of the Casimir functions of $c$ on $M$.
We take $\Span$ in order  to
obtain a vector space: a sum of two Casimir functions for different
$c_1,c_2\in J_0$ need not be a Casimir function.
\end{defi}

The following proposition shows how  the degenerate
bihamiltonian structures can be applied for constructing the completely
integrable systems.

\abz\label{a10.250}
\begin{prop}
Let $J$ be a degenerate bihamiltonian structure on $M$. A family
${\cal F}_0$ is involutive
with respect to any $c_\la\in J$.
\end{prop}

\Proof Let $c_1,c_2\in J_0$ be linearly independent,
$f_i\in Z_{c_i}, i=1,2$. Then
\sabz\label{f1}
\begin{equation}
\{f_1,f_2\}_{c_\la}=(\la_1c_1(f_1)+\la_2c_2(f_1))f_2=-\la_2c_2(f_2)f_1=0.
\end{equation}
Now it remains to prove that for any $c\in J_0,f_i\in Z_c,
i=1,2$, one has
$\{f_1,f_2\}_{c_\la}=0$.
For that purpose we first rewrite
(\ref{f1}) as
\sabz\label{f2}
\begin{equation}
c_\la(x)(\phi_1,\phi_2)=0,
\end{equation}
where $\phi_i\in \ke c_i(x),\ i=1,2,\ x\in M$, and the
lefthandside denotes the contraction of the bivector with two covectors.
Second, we fix $x$ such that $\ra c(x)=R_0$
and approximate $df_2|_x$ by a sequence of elements
$\{\phi^i\}_{i=1}^\infty,\ \phi^i\in \ke c^i(x)$, where $c^i\in
J_0,\ i=1,2,\dots,$ is linearly independent with $c$. Finally,
by (\ref{f2}) we get $c_\la(x)(df_1|_x,\phi^i)=0$ and by the
continuity $\{f_1,f_2\}_{c_\la}(x)=0$. Since the set of such
points $x$ is dense in $M$, the proof is finished. \endproof

In fact this proposition is true for the local Casimir
functions (for the germs of Casimir functions).
The corresponding family of functions (germs)
$\Span_\R (\bigcup_{c\in J_0}Z_c(U))$
($\Span_\R (\bigcup_{c\in J_0}Z_{c,x}$) is denoted by ${\cal F}_0(U)$
(${\cal F}_{0,x}$).

In order to obtain a completely integrable
system from Casimir functions one should
require additional assumptions on the bihamiltonian structure $J$.
Off course, the condition of completeness
given below concerns the local Casimir functions (in fact their germs)
and may be insufficient for
obtaining the completely integrable system. However, it is of
use if the local Casimir functions are restrictions of the
global ones (see Example \ref{a20.50}, below).

Given a Poisson bivector $c_\la\in J$, let $S_\la(x)$ denote the
symplectic leaf of $c_\la$
through a point $x\in M$.

\abz\label{a20.10}
\begin{defi} \rm
(\cite{bols}) Let $J$ be a bihamiltonian structure;
fix some $c_\la\in J$.

$J$ is called {\em complete at a point $x\in M$ with
respect to $c_\la$} if the
linear subspace of $T_x^*M$ generated by the
differentials of
the germs $f\in {\cal F}_{0,x}$ restricted to $S_\la(x)$ has dimension
$\frac{1}{2}\dim S_\la(x)$.
\end{defi}

\abz\label{a20.20}
\begin{prop}
A bihamiltonian structure
$J$ is complete with respect to $c_\la\in J_0$ at a point $x\in
M$ such that $S_\la(x)$ is of maximal dimension
if and only
if $\dim(\bigcap_{c_\la\in J_0}\!
T_x S_\la(x))\lbr =\frac{1}{2}\dim S_\la(x)$.
\end{prop}

The following theorem is due to A.Brailov (see \cite{bols}, Theorem 1.1 and
Remark after it).

\abz\label{a20.30}
\begin{thm}
A bihamiltonian structure
$(J,c_1,c_2)$ is complete with respect to $c_\la\in J_0$ at a point $x\in M$
such that $S_\la(x)$ is of maximal dimension if and only
if the following condition holds
\begin{description}
\item[$(*)$] $\ra (\la_1c_1+\la_2c_2)(x)=R_0$ for any
$\la =(\la_1,\la_2)\in \C^2\setminus \{0\}$.
\end{description}
\end{thm}

Here the bivector $c_\la=(\la_1c_1+\la_2c_2)(x)$ is regarded as an element of
$\bigwedge^2T_x^\C M$, where $T^\C M$ is the complexified tangent bundle, and
its rank is defined as that of the associated sharp map
$c_\la ^\sharp (x):(T_x^\C M)^*\longrightarrow T_x^\C M$.

The theorem shows that $J$ is complete with respect to a fixed
$c_\la\in J_0$ at a point $x$ such that the dimension $S_\la(x)$ is
maximal
if and only if $J=J_0\bigcup\{0\}$ and $J$ is complete at $x$ with
respect to any nontrivial $c_\la\in J$. This motivates the next definition.

\abz\label{a20.40}
\begin{defi} \rm
Let $(J,c_1,c_2)$ be a bihamiltonian structure.
The structure $J$ (the pair $(c_1,c_2)$)
is {\em complete at a point $x\in M$} if condition
$(*)$ of Theorem \ref{a20.30} holds at $x$. $J$ ($(c_1,c_2)$)
is called {\em complete} if it is so at any point from some
open and dense subset in $M$.
\end{defi}

\abz\label{a20.45}
\begin{prop}
Let $J$ be  complete on $M$ and let $x\in M$ be a point of completeness. Then
there exists a neighbourhood $U\ni x$ such that the foliation ${\cal L}$
defined on $U$ by ${\cal F}_0(U)$ is lagrangian in any $S_\la(y),\la\ne
0,y\in U$ (by Proposition \ref{a20.20} this foliation can be defined as the
intersection of the foliations of symplectic leaves for $c_\la\in J_0$).
\end{prop}

\abz\label{a20.48}
\begin{defi}
\rm Call ${\cal L}$ a {\em bilagrangian foliation} of $J$.
\end{defi}

\abz\label{a20.50}
\begin{exa}
\rm (Method of argument translation, see \cite{mf},
\cite{bols}.) Let $\A$ be a Lie algebra, $\D$
its dual space. Fix a basis $\{e_1,\ldots,e_n\}$ in $\A$
with the structure constants $\{c_{ij}^k\}$; write
$\{e^1,\ldots,e^n\}$ for the  dual basis in $\D$.
The standard linear Poisson bivector on $\D$ is defined as
$$c_1(x)=c_{ij}^kx_k\frac{\d}{\d x_i}\wedge\frac{\d}{\d x_j},$$
where $\{x_k\}$ are linear coordinates in $\D$
corresponding to $\{e^1,\ldots,e^n\}$. In more invariant terms
$c_1$ is described as an operator dual to the Lie-multiplication map $[\,,]:
\A\wedge\A\longrightarrow \A$. It is well-known
that the symplectic leaves of $c_1$ are the coadjoint
orbits in $\D$. Now define $c_2$ as a bivector with
constant coefficients $c_2=c(a)$, where $a$ is a fixed point on
any leaf of maximal dimension. It turns out that $c_1,c_2$ form
a Poisson pair and it is easy to describe the set $I$ of points
$x$ for which condition $(*)$ fails. Consider the
complexification $(\D)^\C\cong (\A^{\bf
C})^*$ and the sum $\Sing(\A^\C )^*$ of symplectic
leaves of nonmaximal dimension for the complex linear bivector
$c_{ij}^kz_k\frac{\d}{\d z_i}\wedge\frac{\d}{\d z_j},$
where $z_j=x_j+\i y_j,\ j=1,\ldots,n,$ are the corresponding
complex coordinates in $(\D)^\C$. Then $I$ is equal
to the intersection of the sets $\A^*\subset(\D)^\C$ and
$\overline{a,\Sing(\A^\C )^*}$,
where $\overline{a,\Sing(\A^\C )^*}$ denotes a cone of
complex $2$-dimensional subspaces passing through $a$ and
$\Sing(\A^\C )^*$.

In particular, $(c_1,c_2)$ is complete for a semisimple $\A$
($\codim\Sing(\A^\C)^*\ge 3$, see \cite{adams}, Corollary
4.42,  and codimension of $I$ in $\D$ is not less than $2$).
Note, that this gives rise to completely
integrable systems since the local Casimir functions on $\D$
are restrictions of the global
ones, i.e. the invariants of the coadjoint action.
\end{exa}

\abz\label{a20.60}
\begin{exa}
\rm (Bihamiltonian structure of general position on
an odd-dimensional manifold, see \cite{gz2}.) Consider a pair of
bivectors $(a_1,a_2)$, $a_i\in\bigwedge^2V,  i=1,2$, where $V$ is a
$(2m+1)$-dimensional vector space; $(a_1,a_2)$ is in general
position if and only if is represented by the Kronecker block of
dimension $2m+1$, i.e.
\begin{equation}
\begin{array}{l}
a_1=p_1\wedge q_1+p_2\wedge q_2+\cdots+p_m\wedge q_m\\
a_2=p_1\wedge q_2+p_2\wedge q_3+\cdots+p_m\wedge q_{m+1}
\end{array}
\end{equation}
in an appropriate basis $p_1,\ldots p_m,q_1,\ldots q_{m+1}$ of $V$.
A bihamiltonian structure $J$ on a
$(2m+1)$-dimensional $M$ is in general position if and only if
the pair $(c_1(x),c_2(x))$ is so for any $x\in M$. Such $J$ is
complete. In general, a complete
Poisson pair at a point is a direct sum of the Kronecker
blocks as the corollary of the next theorem shows.
This theorem is a reformulation of the classification result for
pairs of $2$-forms in a vector space
(\cite{gz1}, \cite{gz3}).
\end{exa}

\abz\label{a20.65}
\begin{thm}
Given a finite-dimensional vector
space $V$ over $\C$ and a pair of bivectros $(c_1,c_2),\
c_i\in\bigwedge^2V,$ there exists a direct decomposition
$V=\oplus V_j,\ c_i=\sum c_i^{(j)},
\ c_i^{(j)}\in\bigwedge^2V_j,\ i=1,2,$ such that each triple
$(V_j,c_1^{(j)},c_2^{(j)})$ is from the following list:
\begin{description}
\item[(a)] the Jordan block: $\dim V_j=2n_j$ and in an
 appropriate basis of $V_j$ the matrix of $c_i^{(j)}$ is
equal to $$\left(\begin{array}{cc} 0&A_i\\ -A_i^T&0
\end{array}\right),i=1,2,$$
where $A_1=I_{n_j}$ (the unity $n_j\times n_j$-matrix) and
$A_2=J_{n_j}^\la$ (the Jordan block with the eigenvalue $\la$);
\item[(b)] the Kronecker block: $\dim V_j=2n_j+1$ and in an
appropriate basis of $V_j$ the matrix of $c_i^{(j)}$ is
equal to $$\left(\begin{array}{cc} 0&B_i\\ -B_i^T&0
\end{array}\right), i=1,2,$$
where
$%\renewcommand{\arraycolsep}{-.005cm}
B_1=\left(\begin{array}{cccccc}
1&0&0&\ldots&0&0\\
0&1&0&\ldots&0&0\\
 & & &\ldots& & \\
0&0&0&\ldots&1&0
\end{array} \right),B_2=\left(\begin{array}{cccccc}
0&1&0&\ldots&0&0\\
0&0&1&\ldots&0&0\\
 & & &\ldots& & \\
0&0&0&\ldots&0&1
\end{array}\right)$ ($(n_j+1)\times n_j$-matrices; in case
$n_j=0$ put $c_1=c_2=0$).
\end{description} \end{thm}

\abz\label{a20.67}
\begin{cor}
Let $(J,c_1,c_2)$ be a bihamiltonian structure. It is complete at a point
$x\in M$ if and only if the pair $(c_1(x),c_2(x)),\
c_i(x)\in\bigwedge^2(T_x^{\C}M), i=1,2,$ does not contain the
Jordan blocks in its decomposition.
\end{cor}

\Proof The statement follows from the definition of completeness. \endproof

\abz\label{20.68}
\begin{rem}
\rm
In the absence of the Jordan blocks Theorem \ref{a20.65} is valid also over
reals. 
\end{rem}

\section{Complete bihamiltonian structure at a point.}
Now, we shall examine a linear bihamiltonian structure $(J,c_1,
c_2), c_i\in \bigwedge^2 V$, where $V$ is a vector space over $\R$,
such that the decomposition $V=\oplus_{j=1}^kV_j,
\ c_i=\sum_{j=1}^kc_i^{(j)}$
(see Theorem \ref{a20.65} and Remark \ref{20.68}) consists of
$k$ Kronecker blocks $V_1,\ldots,V_k,\lbr \dim V_j=2n_j+1,
n_1<\ldots < n_k$.

The aim is to extract the invariants and to introduce the
infinithesimal approximation to Veronese webs (these last
will be defined in the next section).

It turns out that the decomposition to Kronecker blocks is
noninvariant. To illustrate this let us
consider $V=\Span \{ e,p,q_1,q_2\} ,c_1=p\wedge q_1,
c_2=p\wedge q_2$. Here $V=V_1\oplus V_2$, where $V_1=\Span
\{ e\}, V_2=\Span\{ p,q_1,q_2\}$, but instead $V_1$ one can
choose any direct complement to $V_2$. However, there is a
canonically defined filtration associated to $J$.

Let $P_{c_\la}\subset V$ ($P_{c_\la^{(j)}}\subset V_j$) be
the characteristic subspace, i.e. the symplectic leaf
through $0$, of $c_\la =\la _1c_1+\la _2c_2$
($c^{(j)}_\la =\la _1c^{(j)}_1+\la _2c^{(j)}_2$),
$(\la _1,\la _2)\in\R ^2$; let $L=\cap _{\la\ne 0}
P_{c_\la}$ ($L_j=\cap _{\la\ne 0}P_{c^{(j)}_\la}$) be the
bilagrangian subspase, i.e. the leaf through $0$ of the
bilagrangian foliation, corresponding to the
bihamiltonian structure $J$ ($\{c^{(j)}_\la\}_{\la\in\R
^2}$), see Definition \ref{a20.48}.

Put
$$\Phi_i=\sum_{\mbox{distinct}\ \la_1,\ldots,\la_i\in\pr (\R
^2)}P_{c_{\la_1}}\cap\ldots \cap P_{c_{\la_i}},
i=1,2,\ldots , \Phi _0=V.
$$

\abz\label{30.10}
\begin{thm}
The following relations hold
$$
\Phi_0=\Phi_1=\cdots=\Phi_{n_1}\supset\Phi_{n_1+1}=\Phi_{n_1+2}=\cdots=
\Phi_{n_2}\supset\cdots\supset\Phi_{n_{k-1}+1}=\cdots=\Phi_{n_k}\supset
$$
$$
\supset\Phi_{n_k+1}=\Phi_{n_k+2}=\ldots=:\Phi_{n_{k+1}},
$$
where
$$
\Phi_{n_j}=\sum_{ l< j}L_l\oplus\sum_{l\ge j}V_l,
j=1,\ldots ,k+1
$$ (we put $V_l=0, l>k$). In particular,
the filtration is stabilized from $i>n_k$,
$\Phi_{n_k+1}=\oplus L_j=L$, and the numbers
$n_1,\ldots,n_k$ are invariants of $J$.
\end{thm}

\abz\label{30.15}
\begin{rem} \rm
The filtration $F_0=\Phi_{n_1}^\bot\subset\cdots\subset
F_{k-1}=\Phi_{n_k}^\bot\subset F_{k}=L^\bot=(V/L)^*$ ($\bot$
stands for the annihilator sign) appears in \cite{gz5} and
is called there isotypic.  We shall refer to this notion
below.
\end{rem}

Before we begin to prove the theorem we recall the following definition.

\abz\label{30.30} \begin{defi} \rm
(\cite{gz2}) Let ${\cal S},V$ be vector spaces of
dimensions $2$    and $n+1, n\ge 0$, respectively. A {\em
Veronese inclusion} of $\pr ({\cal S})$ in $\pr (V)$ is a
map $i:\pr ({\cal S})\longrightarrow\pr (V)$ such that
there exists a linear isomorphism $\phi:\pr
(V)\longrightarrow\pr (S^n{\cal S})$ making the following
diagram commutative:

$$
\begin{array}{ccc}
\pr ({\cal S})  & \stackrel{i}{\longrightarrow}& \pr (V) \\
\parallel &  & \uparrow  \phi \\
\pr ({\cal S}) &\stackrel{ \pr (S^n(\cdot))}{\longrightarrow}
&\pr (S^n{\cal S}).
\end{array}
$$
The image $i(\pr ({\cal S}))$ is called {\em Veronese
curve}.
\end{defi}

Here $S^n$ denotes the $n$-th symmetric power; the standard model of
the mapping $S^n(\cdot)$ is described as follows. Let ${\cal S}$
be a space of linear functions $f$ in two variables $t_1,t_2$.
Then $S^n{\cal S}$ is a space of homogeneous polynomials in
$t_1,t_2$ and $S^n(f)=f^n$. For $n=0$ the map $i$ is not an inclusion, but
we shall use the defined term in this situation as well.

\Proof We now shall prove Theorem \ref{30.10}. It is sufficient to show the
following equalities
\sabz\label{form}
\begin{equation}
\Phi_iV_j:=\sum_{\mbox{distinct}\ \la_1,\ldots,\la_i}P_{c^{(j)}_{\la_1}}
\cap\ldots\cap P_{c^{(j)}_{\la_i}}=
\left\{\begin{array}{ll}
L_j&i\ge n_j+1\\
V_j&i<n_j+1.
\end{array}
\right.
\end{equation}
One has
$$
(\Phi_iV_j)^\bot=\bigcap_{\mbox{distinct}\ \la_1,\ldots,\la_i}
(P^\bot _{c^{(j)}_{\la_1}}
+\cdots + P^\bot _{c^{(j)}_{\la_i}}).
$$
The 1-dimentional annihilator $P^\bot_{c^{(j)}_\la}\in \pr
((V_j/L_j)^*)$ sweeps an appropriate Veronese curve (see \cite{gz2}).
On the other hand, images of distinct points under a Veronese inclusion $\pr
({\cal S})\rightarrow\pr  (V)$ are
linearly independent untill their number does not exceed
$\dim V$ (cf.~\cite{cg2}, I.A).  So now, the first of
equalities \ref{form} follows from the fact that $P^\bot
_{c^{(j)}_{\la_1}} +\cdots + P^\bot
_{c^{(j)}_{\la_i}}=(V_j/L_j)^*$ for any set of
distinct $\la_1,\ldots,\la_i, i\ge n_j+1$.

For the second one, we notice that for any set of points $\la_1,\ldots
,\la_{n_j+1}$
$$
\bigcap_{s=1}^{n_j+1}(P_1
+\cdots + \hat{P_s}+\cdots +P_{n_j+1})=\{ 0\},
$$
where we put $P_s=P^\bot _{c_{\la_s}}$ and $\hat{}$ means
omitting of the corresponding term.  Thus we
proved it for $i=n_j$; this implies \ref{form} also for $i<
n_j$.  \endproof

\abz\label{30.40}
\begin{cor}
Let $0=F_1\subset \cdots\subset
F_{k-1}\subset F_{k}=(V/L)^*$ be the isotypic filtration (see \ref{30.15}).
Put $F_{j}P^\bot_{c_\la}=F_{j}
\cap P^\bot_{c_\la},j=1,\ldots,k$. Then:
1) $A_{j}^\la:=F_{j}P^\bot_{c_\la}/F_{j-1}P^\bot_{c_\la}$
is a one-dimensional subspace in $A_j:=F_{j}/F_{j-1}$; 2)
the mapping $\pr (\R^2)\ni\la\mapsto A_j^\la\in\pr (A_j)$
is a Veronese inclusion for any $j=1,\ldots,k$.  \end{cor}

\Proof We first notice that
$F_{j}=\oplus_{i=1}^{j}(V_i/L_i)^*,j=1,\ldots,k$, as
Theorem \ref{30.10} implies. Under the identification
$(V/L)^*= \oplus_{i=1}^k(V_i/L_i)^*$
$$
F_{j}P_{c_\la}^\bot=\oplus_{i=1}^jP_{c_\la^{(i)}}^{\bot_i},
$$
where $\bot_i$ stands for the annihilator of a subspace in $V_i/L_i$. Thus
there are linear isomorphisms
$A_j\cong(V_j/L_j)^*$ and $A_j^\la\cong P_{c_\la^{(j)}}^{\bot_j},
j=1,\ldots,k$, and 2) follows from the analogous fact for a
sole Kronecker block (\cite{gz2}).  \endproof

\abz\label{rema}
\begin{rem}
\rm We note that analogues of Theorem~\ref{30.10} and
Corollary~\ref{30.40} can be proved also without the
restriction that $n_1,\ldots,n_k$ are distinct. In order to
do that one should use a "multiple" version of a Veronese
inclusion, i.e. a map $\phi :\pr ({\cal S})\rightarrow
G(k,V^{(l+1)k})$ ($G(k,V)$ denotes the Grassmannian of
$k$-planes in a $(l+1)k$-dimensional vector space $V$) such
that there exist a decomposition $V=\oplus_{j=1}^kV_j,
\dim V_j=l+1$, and Veronese inclusions $i_j:\pr({\cal
S})\rightarrow\pr(V_j),j=1,\ldots,k$, with the property
$\phi(v)=\Span_\R\{ i_1(v),\ldots,i_k(v)\}$, where $i_j(v)$
is considered as a 1-dimensional subspace in $V$. This
definition can be also used to adapt the notion of a
Veronese web (see Definition~\ref{40.20}, below) to a more
general situation.
\end{rem}

\abz\label{30.50}
\begin{defi} \rm
An {\em infinithesimal Veronese web} of type
$(n_1,\ldots,n_k),n_1<\ldots<n_k$,  on a vector
space $W, \dim W=n_1+\cdots+n_k+k$, is a 1-parameter family
$\{\W_\la\}_{\la\in\pr(\S)}$
of linear subspaces $\W_\la\subset W, \codim\W_\la=k$, satisfying the
following conditions:
\begin{description}
\item[(i)] there is a filtration $0=F_{0}\subset\ldots\subset
F_{k-1}\subset F_{k}=W^*$ of the dual space with $\dim
F_{j}/F_{j-1}=n_j+1, j=1,\ldots,k$;
\item[(ii)] it induces the filtration $0=F_{0}{\cal
W}^\bot_\la\subset\ldots\subset F_{k-1}{\cal
W}^\bot_\la\subset F_{k}{\cal
W}^\bot_\la ={\cal
W}^\bot_\la, F_{j}{\cal
W}^\bot_\la =F_{j}\cap{\cal
W}^\bot_\la , j=1,\ldots,k$, of the annihilator ${\cal
W}^\bot_\la\subset W^*$ so that $\dim F_{j}{\cal
W}^\bot_\la =j$; in particular $A_j^\la:=F_{j}{\cal
W}^\bot_\la/F_{j-1}{\cal
W}^\bot_\la$ can be considered as 1-dimensional subspace in $A_j:=
F_{j}
/F_{j-1}$;
\item[(iii)] the map $\pr({\cal S})\ni\la\mapsto A_j^\la\in\pr(A_j)$
is a Veronese inclusion, $j=1,\ldots,k$.
\end{description}
\end{defi}

\abz\label{30.60}
\begin{prop}
Let $(J,c_1,c_2)$ be as above. Then the vector space
$W=V/L$ has a structure of an infinithesimal Veronese web of
type $(n_1,\ldots,n_k)$.
\end{prop}

\Proof The proof follows from Corollary \ref{30.40}.
\endproof

\section{Simple bihamiltonian structures and
their Veronese webs.}
In this section we shall define objects that generalize the
Veronese webs introduced in \cite{gz2} for the bihamiltonian structures of
general position. We shall show that any complete bihamiltonian structure
from the  class defined below
admits the local reduction to such an object.

\abz\label{40.10}
\begin{defi} \rm
Let $J$ be a complete bihamiltonian structure on $M$.
A {\em type} of $J$ at $x\in M$ is the vector
$(n_1,\ldots,n_k)(x)$, where
$2n_1(x)+1,\ldots,2n_k(x)+1$ are
dimensions of the Kronecker blocks in the decomposition of $(c_1(x),c_2(x))$
for some generating $J$ Poisson pair $(c_1,c_2)$ (these dimensions do not
depend on this pair, see Theorem \ref{30.10}).  If this vector is
independent of $x$ we call it a type of $J$ and say that $J$ is {\em
regular} (cf. Example~\ref{ex}, below). If, moreover, all
$n_j, j=1,\ldots,k$, are different we call $J$ {\em
simple}.  \end{defi}

\abz\label{40.20}
\begin{defi} \rm
Consider a manifold $U$ diffeomorphic to an open set in $\R^N$, where
$N=(n_1+1)+\cdots+(n_k+1),
n_1<\ldots<n_k$, and a family
$\W=\{\W_\la\}_{\la\in\pr({\cal S})}$ of $k$-codimensional
foliations on $U$ parametrized by the projectivizaton of a
two-dimensional vector space ${\cal S}$.  We call $\W$ a
{\em Veronese web} of type $(n_1,\ldots,n_k)$ if the following  conditions are
satisfied:
\begin{description} \item[(i)] there is a
bundle filtration $0=F_{0}\subset \ldots\subset F_{k-1}\subset
F_{k}=T^*U$ such that  \linebreak $\ra F_{j}/
F_{j-1}=n_{j}+1, j=1,\ldots,k$; 
\item[(ii)] it induces
the bundle filtration $0=F_{0}{\cal
W}_\la^\bot\subset\ldots\subset F_{k-1}{\cal
W}_\la^\bot\subset F_{k}{\cal W}_\la^\bot ={\cal
W}_\la^\bot, \lbr F_{j}{\cal W}_\la^\bot=F_{j}\cap {\cal W}_\la^\bot,
j=1,\ldots,k$, of the annihilating bundle ${\cal W}_\la^\bot:=(T{\cal
W}_\la)^\bot\subset T^*U$ so that $\ra F_{j}{\cal W}_\la^\bot=j$; in
particular $A_j^\la(x):=\lbr F_{j}{\cal W}_{\la,x}^\bot /
F_{j-1}{\cal W}_{\la,x}^\bot$ can be considered as a
1-dimensional subspace in $A_j(x):= F_{j,x}/ F_{j-1,x}$
for any $x\in U$; 
\item[(iii)] the map $\pr({\cal
S})\ni\la\mapsto A_j^\la(x)\in\pr(A_j(x))$ is a Veronese
inclusion for any $x\in U, j=1,\ldots,k$.
\end{description}
\end{defi}

\abz\label{40.30}
\begin{thm}
Let $J$ be a simple bihamiltonian structure of type ${\bf
n}=(n_1,\ldots,n_k),\lbr  n_1<\cdots<n_k$, and let $x\in M$ be a point of
completeness for $J$. Write $\V_\la$ for the foliation of
symplectic leaves of $c_\la\in J$. Then there exists a
neighbourhood $\tilde{U}\ni x$ such that $U=\tilde{U}/{\cal L}$
(see \ref{a20.45}) is diffeommorphic to an open set in $\R^N$
and $\{\V_\la|_{\tilde{U}}/{\cal L}\}_{\la \in\pr(\S)}$ is a
Veronese web of type ${\bf n}$ on $U$.
\end{thm}

\Proof The theorem follows from Proposition \ref{30.60}.\endproof

\abz\label{40.40}
\begin{exa}
\rm  Let $U=\R^3(x,y,z)$, $\alpha _1=xdy-dz, \alpha_2^\la=\la_1dx+\la_2dy,
k=2,n_1=0,n_2=1$.  Put $\G(F_{1})=\Span\{\alpha_1\}, \G((T{\cal
W}_\la)^\bot)=\Span\{\alpha_1,\alpha_2^\la\}$, where $\G$ stands for the
space of sections and $\Span$ is taken over the ring of functions. Then
$\G(F_{1}(T{\cal W}_\la)^\bot)=\G(F_{1})$. Since $T{\cal W}_\la\subset
TU$ is a subbundle of rank 1, it is indeed tangent to
1-dimensional foliation ${\cal W}_\la$. Explicitely,
$\G(T{\cal W}_\la)=\Span\{\la_2v_1-\la_1v_2\}$, where
$v_1=\frac{\d}{\d x}, v_2=\frac{\d}{\d y}+x\frac{\d}{\d
z}$.  On $\tilde{U}=\R (p)\times U$ one defines the
corresponding bihamiltonian structure  as $\{\frac{\d}{\d
p}\wedge(\la_1v_1+\la_2v_2)\}_{(\la_1,\la_2)\in\R^2}$.
\end{exa}

\abz\label{40.50}
\begin{exa}
\rm Let $U=\R^5(x,y,z,s,t)$,
$\alpha_1^\la=\la_1(xdy-dz)+\la_2(sdy-dt),
\alpha_2^\la=\la_1^2dx+\la_1\la_2ds+\la_2^2dy,
\G(F_{1})=\Span\{xdy-ds,sdy-dt\}, \G((T{\cal
W}_\la)^\bot))=\Span\{\alpha_1^\la,\alpha_2^\la\}, k=2,n_1=1,n_2=2$. Then
$\G(F_{1}(T{\cal
W}_\la)^\bot))=\Span\{\alpha_1^\la\}$. The 3-distribution in $TU$
annihilating by 1-forms $\alpha_1^\la,\alpha_2^\la$ is integrable since
$$
d\alpha_1^\la=\left\{\begin{array}{ll}
\alpha_2^\la\wedge\frac{1}{\la_1}dy & \mbox{if}\ \la_1\not
=0\\ -\alpha_2^\la\wedge\frac{1}{\la_2}ds & \mbox{if}\
\la_1=0.  \end{array} \right.  $$ Explicitely, $\G(T{\cal
W}_\la)=\Span\{\la_2\frac{\d}{\d x}-\la_1 \frac{\d}{\d
s},\la_2\frac{\d}{\d s}-\la_1v,\la_2\frac{\d}{\d z}-\la_1
\frac{\d}{\d t}\}$, where $v=\frac{\d}{\d y}+x\frac{\d}{\d z}+s \frac{\d}{\d
t}$, and on
$\tilde{U}=\R^3(p_1,p_2,p_3)\times U$
the corresponding
bihamiltonian structure is
$\{\frac{\d}{\d p_1} \wedge(\la_2\frac{\d}{\d x}-\la_1\frac{\d}{\d
s})+\frac{\d}{\d p_2}\wedge(\la_2\frac{\d}{\d s}-\la_1v)+\frac{\d}{\d
p_3}\wedge(\la_2\frac{\d}{\d z}-\la_1\frac{\d}{\d
t})\}_{(\la_1,\la_2)\in\R^2}$.
\end{exa}

\abz\label{40.55}
\begin{rem}
\rm
Off course, it is more convenient to describe Veronese webs
in terms of a bundle direct decomposition
$B_{1}\oplus\ldots\oplus  B_{k}=T^*U$ such that
$F_{j}=\oplus_{i=1}^jB_{i}$ rather than in terms of the
isotypic filtration itself. In the above examples we used implicitely
such a decomposition. However, one should remember that it
is not unique. For instance, in Example \ref{40.50} one has
$B_{1}=F_{1}, \G(B_{2})=\Span\{ dx,ds,dy\}$. But one could
take
$\tilde{\alpha}_2^\la=\la_1^2(dx+xdy-dz)+\la_1\la_2(ds+sdy-dt)
+\la_2^2dy=\alpha_2^\la+\la_1\alpha_1^\la$ instead of $\alpha_2^\la$
and $\G(\tilde{B}_{2})=\Span\{ dx+xdy-dz,ds+sdy-dt,dy\}$.
Although this does not change the web, the corresponding
decomposition is changed. In \cite{gz5} the author gives an involved  
analysis of this
nonuniqueness.   \end{rem}

\abz\label{40.60}
\begin{defi} \rm
A Veronese web admits the following local description. One can
choose linear coordinates $(\la_1,\la_2)$ on $\S$ and a local
coframe $\alpha_1^1,\ldots,\alpha_{n_1+1}^1,\ldots,
\alpha_1^k,\ldots,\lbr \alpha_{n_k+1}^k$ such
that
$\alpha_1^j,\ldots,\alpha_{n_j+1}^j\in\G(F_{j}),j=1,\ldots,k$, and
the annihilator $(T\W_\la)^\bot\subset T^*U$ is
generated by $\alpha_\la^1,\ldots,\alpha_\la^k$, where
$\alpha_\la^j=\la_1^{n_j}\alpha_1^j+
\la_1^{n_j-1}\la_2\alpha_2^j+\cdots
+\la_2^{n_j}\alpha_{n_k+1}^j$ (Veronese curve). If in a
neighbourhood of any $x\in U$ there exists a holonomic coframe
with the above properties, the Veronese web is called
{\em flat}.

In particular, all bundles in the isotypic filtration of a flat web are
completely integrable as differential systems and, moreover, such a web
splits to a direct product of flat Veronese webs of codimension $1$.
\end{defi}

The webs from Examples \ref{40.40}, \ref{40.50}  are not flat, since the
bundles $F_{1}$ are nonintegrable.

We conclude the section by an example of a complete
bihamiltonian structure that is not regular.

\abz\label{ex}
\begin{exa}
\rm  Let $M=\R^6$ with coordinates
$(p_1,p_2,q_1,\ldots,q_4)$, $c_1=\frac{\d}{\d
p_1}\wedge\frac{\d}{\d q_1}+\frac{\d}{\d
p_2}\wedge\frac{\d}{\d q_2}, c_2=\frac{\d}{\d
p_1}\wedge(\frac{\d}{\d q_2}+q_1\frac{\d}{\d q_3})+\frac{\d}{\d
p_2}\wedge\frac{\d}{\d q_4}$. Here we have: two $3$-dimensional
Kronecker blocks on $M\setminus H, H=\{q_1=0\}$; the
$5$-dimensional Kronecker block and the $1$-dimensional zero
block on the hyperplane $H$.
\end{exa}

\section{Veronese webs for the argument translation method.}
The notations from  Subsection \ref{a20.50} will be used below.
We consider
normal (d\'{e}ployable in terminology of Bourbaki, \cite{bourbIX}, IX,3)
real form $\A$ of complex simple Lie algebra. Let
$m_1,\ldots,m_r,r=\ra(\A)$ be the exponents of $\A$.

\abz\label{50.10}
\begin{thm}
Let $(c_1,c_2)$ be the Poisson pair from Example \ref{a20.50}.
Then the Veronese web
$\{\W_\la\}_{\la\in\R^2}$ of the
corresponding bihamiltonian structure $J$
is of type $(m_1,\ldots ,m_r)$ and is flat (Definition
\ref{40.60}) in a neighbourhood of any point $\pi (x)$,
where $x\in(\D\setminus I)$ and $\pi$ denotes the canonical
projection $\pi :\D\setminus I\rightarrow
(\D\setminus I)/{\cal L}$ (cf. \ref{a20.45}, \ref{40.30}).
\end{thm}

\Proof   Let
$g_1(x),\ldots,g_r(x), \deg g_j=m_j+1$, be a set of
algebraically independent  global homogeneous
polynomial Casimir functions for $c$  (see \cite{bourb}, VIII,8).
Here we have identified $\A$ and $\D$ by means of the Killing form.
Note that $g_1,\ldots,g_r$ are functionally
independent on $\A\setminus\Sing\, \A$, where $\Sing\,\A$ is the
set of adjoint orbits of nonmaximal dimension.
Indeed, their restrictions to a Cartan subalgebra ${\frak h}\subset\A$
are algebraically independent and invariant with respect to the
Weyl group $W$. Now, we can apply the
result of R.Steinberg (\cite{steinberg}) to deduce the
nondegeneracy for the Jacobi matrix of $g_1|_{\frak
h},\ldots ,g_r|_{\frak h}$ at a regular point.

Consider the subspace $d{\cal F}_0\subset
\G(T^*\D)$ generated by the differentials of functions from
the involutive set ${\cal F}_0$ (see \ref{a10.230})
corresponding to $J$. It turns out
that $d{\cal F}_0$ is generated by
$\{dg_j|_{\la_1x+\la_2a},(\la_1,\la_2)\in\R^2,j=1,\ldots,r\}$.
If $g_j^i(a,x),i=0,\ldots,m_j+1,j=1,\ldots,r$, are the coefficients
of the  Taylor expansions $g_j(x+\la
a),j=1,\ldots,r$, with respect to $\la \in\R$, then one also has
\sabz\label{f}
\begin{equation}
d{\cal F}_0=\Span\{dg_j^i(a,x),i=0,\ldots,m_j,j=1,\ldots,r\}.
\end{equation}

Moreover, these differentials  are linearly independent at any
$x\in \D\setminus I$. This follows from the fact that $J$ is
complete at $\D\setminus I$, from (\ref{f}), and from the formula
$\sum_{j=1}^rm_j=\frac{1}{2} (\dim\A -r)$ (cf. \cite{vo}, formula
(F1), p. 289).

Thus, we can regard $g_j^i(a,x),i=0,\ldots,m_j,j=1,\ldots,r$ as
coordinates on the reduced space $(\D\setminus I)/{\cal
L}$.  Finally, $(T\W_\la)^\bot,\la=(\la_1,\la_2)$, is
generated by $$
\la_1^{m_j}dg_j^0(a,x)+\la_1^{m_j-1}\la_2dg_j^1(a,x)+\cdots
+\la_2^{m_j}dg_j^{m_j}(a,x), j=1,\ldots,r.
$$
\endproof

{\nine

}
\end{document}